\documentclass[11pt, twoside]{article}
\usepackage{float}
\usepackage[centertags]{amsmath}
\usepackage{newlfont}
\usepackage{array}
\usepackage{amsmath, amsthm, amscd, amsfonts, amssymb, graphicx, color, delarray}
\usepackage[bookmarksnumbered, plainpages, backref]{hyperref}

\date{}

\begin{document}

\centerline{\Large{\bf Some contributions to Regular Polygons}}

\centerline{}

\centerline{\bf {Deniz \"{O}ncel and Murat Kiri\c{s}ci}}

\centerline{}

\centerline{Department of Mathematical Education, Hasan Ali Y\"{u}cel Education Faculty,}

\centerline{ Istanbul University, Vefa, 34470, Fatih, Istanbul, Turkey}

\centerline{}

\begin{abstract}
The aim of this work is to use Napoleon's Theorem in different regular polygons,
and decide whether we can prove Napoleon's Theorem is only limited with triangles
or it could be done in other regular polygons that can create regular polygons.
\end{abstract}

\centerline{}



\theoremstyle{plain}
\newtheorem{thm}{Theorem}[section]
\numberwithin{equation}{section}
\numberwithin{figure}{section}  
\theoremstyle{plain}
\newtheorem{pr}[thm]{Proposition}
\theoremstyle{plain}
\newtheorem{exmp}[thm]{Example}
\theoremstyle{plain}
\newtheorem{cor}[thm]{Corollary} 
\theoremstyle{plain}
\newtheorem{defin}[thm]{Definition}
\theoremstyle{plain}
\newtheorem{lem}[thm]{Lemma} 
\theoremstyle{plain}
\newtheorem{rem}[thm]{Remark}
\numberwithin{equation}{section}

\section{Introduction}

The famous theorem of Napoleon is one of the most interesting assertions from elementary geometry of planar figures.
Although over 150 years have passed, mathematicians still continues to be of interest to Napoleon's Theorem as both professional and amateur alike.
Because, this theorem is an important and useful tool to expand the mathematical horizon for many mathematicians. In the Euclidean Plane, Napoleon's Theorem is easily proven. A wide variety of proofs of this theorem have been given in a lot of manuscripts.\\

In this work, we give a well known simple proof of Napoleon's Theorem based on Sine and Cosine Laws.
Further, we apply the Theorem to regular polygon such as regular hexagon, square and octagon.\\

Napoleon's Theorem states the fact that if equilateral triangles are drawn outside of any triangle, the centers of the equilateral triangles will form an equilateral triangle.\\

We take a main triangle $\Delta ABC$, three equilateral triangles, which are outside $\Delta ABC$ are $\Delta ADB$, $\Delta BCE$, $\Delta AFC$. We denote the discovered triangle by $\Delta PQR$. Then, $\widehat{FAC}=\widehat{ACF}=\widehat{CFA}=\widehat{ADB}=\widehat{BAD}=\widehat{DBA}=\widehat{BCE}=\widehat{ECB}=\widehat{CBE}=60^{\circ}$ (Figure \ref{Fig. 2}).\\

We give a well known simple proof of Napoleon's Theorem as follows:
\begin{figure}

  \includegraphics[width=70mm]{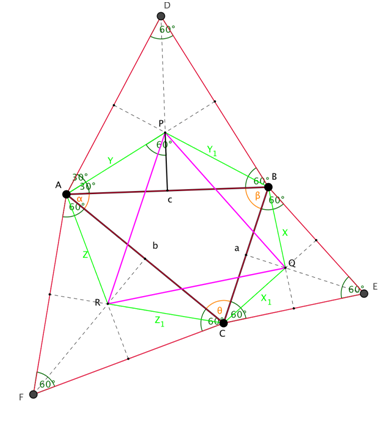}\\
  \caption{}\label{Fig. 2}

\end{figure}

\begin{proof}\cite{Brodie}
Since $\Delta ECB$, $\Delta FCA$, $\Delta ADB$ are equilateral, then
$\widehat{FAC}=\widehat{ACF}=\widehat{CFA}=\widehat{ADB}=\widehat{BAD}=\widehat{DBA}=\widehat{BCE}=\widehat{ECB}=\widehat{CBE}=60^{\circ}$. We take $|CQ|$, $|BQ|$,
$|BP|$, $|AP|$, $|AR|$, $|CR|$ is an angle bisector. Thus, $\widehat{PAB}=\widehat{PBA}=\widehat{QBC}=\widehat{QCB}=\widehat{RCA}=\widehat{RAC}=30^{\circ}$. Isosceles triangles have two equal sides, that is,
$x=x_{1}$, $y=y_{1}$, $z=z_{1}$ and so use the fact that the centroid of an equilateral triangle  $\Delta ADB$, $\Delta BCE$, $\Delta AFC$ lies along each median, $2/3$ of the distance from the vertex to the midpoint of the opposite side,
$x=x_{1}=\frac{2\frac{a\sqrt{3}}{2}}{3}=\frac{a\sqrt{3}}{3}$, $y=y_{1}=\frac{c}{\sqrt{3}}$,
$z=z_{1}=\frac{b}{\sqrt{3}}$. If we use the Cosine formula for $\Delta PQR$,
then, we have

\begin{eqnarray*}
&&|RP|^{2}=y^{2}+z^{2}-2yz\cos(\alpha+60^{\circ})\\
&&|RP|^{2}=\left(\frac{c}{\sqrt{3}}\right)+\left(\frac{b}{\sqrt{3}}\right)-2\frac{c}{\sqrt{3}}\frac{b}{\sqrt{3}}\cos(\alpha+60^{\circ})\\
&&3|RP|^{2}=b^{2}+c^{2}-2bc\cos(\alpha+60^{\circ})\\
&&3|RP|^{2}=b^{2}+c^{2}-2bc(\cos\alpha\cos60^{\circ}-\sin\alpha\sin60^{\circ})\\
&&3|RP|^{2}=b^{2}+c^{2}-2bc(\frac{1}{2}\cos\alpha-\frac{\sqrt{3}}{2}\sin\alpha).
\end{eqnarray*}
If we apply the Law of Cosine to $\Delta ABC$ and the Law of Sine for the area of $\Delta ABC$, then, we can write $a^{2}=b^{2}+c^{2}-2bc\cos\alpha$ and $2(\textit{Area of} \quad \Delta ABC)=bc\sin \alpha$, respectively. Substituting these statement into the $3|RP|^{2}=b^{2}+c^{2}-2bc(\frac{1}{2}\cos\alpha-\frac{\sqrt{3}}{2}\sin\alpha)$ gives
\begin{eqnarray*}
3|RP|^{2}=\frac{1}{2}\left(a^{2}+b^{2}+c^{2}\right)+2\sqrt{3}(\textit{Area of} \quad \Delta ABC).
\end{eqnarray*}

In the same idea, we can compute the length of $|PQ|$ and $|RQ|$. Since, the $\Delta APB$, $\Delta BQC$, $\Delta CRA$ are isosceles triangles, we say $\alpha=\beta=\theta$ and also $|RQ|=|PQ|=|RP|$. Then, it follows that the $\Delta PQR$ connecting the three centroids is equilateral.

\end{proof}

\section{Main Results}

\begin{figure}

  \includegraphics[width=120mm]{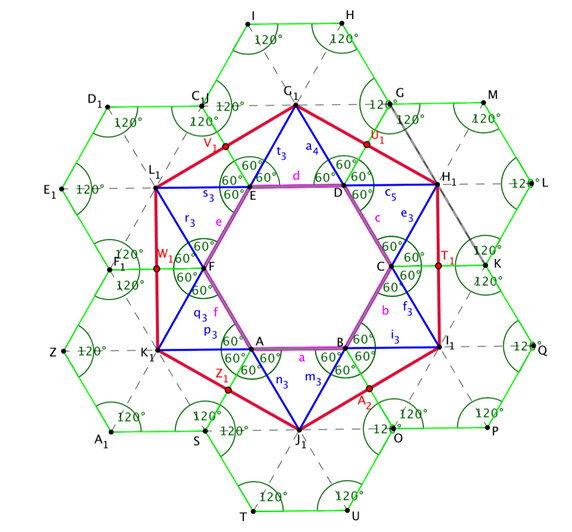}\\
  \caption{}\label{Fig. 3}

\end{figure}

In this section, we will apply Napoleon's Theorem to some regular polygon and prove as in Section 1.\\

\subsection{Application to Regular Hexagon}

Consider the regular hexagon $ABCDEF$ and the six regular hexagons, which are outside the main regular hexagon(Figure \ref{Fig. 3}).\\

Each interior of this hexagon is $120^{\circ}$, because the sum of the interior angles of any hexagon is $720^{\circ}$. Now, we divide
regular hexagons into six pieces by using its diagonals. It will make six equilateral triangles, which are all $60^{\circ}$. Then,
$\widehat{G_{1}ED}=\widehat{G_{1}DE}=<\widehat{H_{1}DC}=\widehat{H_{1}CD}=\widehat{I_{1}CB}=\widehat{I_{1}BC}=\widehat{J_{1}BA}=\widehat{J_{1}AB}=\widehat{K_{1}AF}=\widehat{K_{1}FA}=\widehat{L_{1}EF}=\widehat{L_{1}FE}=60^{\circ}$.\\

The diagonals of regular hexagon are equal. Hence half of the diagonals are equal too. Hence, $c_{5}=e_{3}$, $m_{3}=n_{3}$, $f_{3}=i_{3}$, $t_{3}=a_{4}$, $r_{3}=s_{3}$, $q_{3}=p_{3}$. Also, $a=b=c=d=e$, since regular hexagon has equal sides. We know that the equilateral triangles have equal sides. From this, we can write $c_{5}=e_{3}=c$, $f_{3}=i_{3}=b$, $m_{3}=n_{3}=a$,
$q_{3}=p_{3}=f$, $r_{3}=s_{3}=e$ and $t_{3}=a_{4}=d$. In that case, we obtain $\Delta CH_{1}I_{1}=\Delta BI_{1}J_{1}=\Delta AK_{1}J_{1}=\Delta FL_{1}K_{1}=\Delta EG_{1}L_{1}=\Delta DH_{1}G_{1}$ and $|H_{1}I_{1}|=|I_{1}J_{1}|=|J_{1}K_{1}|=|K_{1}L_{1}|=|L_{1}G_{1}|=|G_{1}H_{1}|$, as we desired.\\

The ratio between the lengths of edges of the main hexagon and the discovered/new hexagon:\\

Let's think that one side of the main hexagon is $x$, $(x\in \mathbb{Z})$.
\begin{eqnarray*}
&&\frac{\sin60^{\circ}}{x}=\frac{\sin60^{\circ}}{|CH_{1}|}\\
&&x=|CH_{1}|\\
&&\Delta H_{1}DC \textit{is isosceles triangle since} |CH_{1}|=|DH_{1}|.
\end{eqnarray*}
Also, $\Delta H_{1}DC$ and $\Delta I_{1}CB$ are equal hence $|DH_{1}|=|BI_{1}|$, the edges $c=c_{5}=e_{3}=x$.
Hence, from the special triangle $(30-60-90)$ in the $\Delta CH_{1}T_{1}$
\begin{eqnarray*}
&&\frac{\sin90^{\circ}}{x}=\frac{\sin60^{\circ}}{|H_{1}T_{1}|}\\
&&\frac{\sqrt{3}}{2}x=|H_{1}T_{1}|.
\end{eqnarray*}
Since $|H_{1}I_{1}|$ is the double of $|H_{1}T_{1}|$ because the edges see the same degree, $60^{\circ}$
\begin{eqnarray*}
|H_{1}I_{1}|=\sqrt{3}x.
\end{eqnarray*}
The one side of the main hexagon is $x$ then, one side of the discovered hexagon need to be $\sqrt{3}x$, since
the ratio between sides won't change no matter the size of the hexagon because same sine law will be applied with the same, equal angles.

\subsection{Application to Square}

\begin{figure}

  \includegraphics[width=120mm]{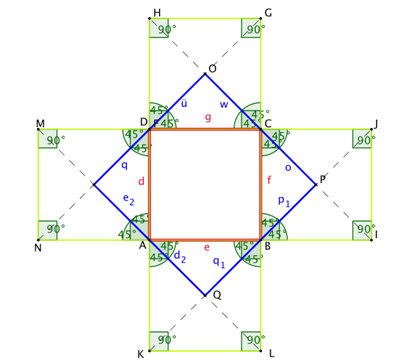}\\
  \caption{}\label{Fig. 4}

\end{figure}

The main square is $ABCD$ and the four squares, which are outside the main square(Figure \ref{Fig. 4}).\\

Each interior of a square is $90^{\circ}$, because the sum of the interior angles of any quadrilateral is $360^{\circ}$.
$DG$ and $CH$, $BJ$ and $CI$, $BK$ and $AL$, $MA$ and $ND$ are the bisector of the squares. Then,
$\widehat{ODC}=\widehat{OCD}=<\widehat{PCB}=\widehat{PBC}=\widehat{QBA}=\widehat{QAB}=\widehat{H_{1}AD}=\widehat{H_{1}DA}=\widehat{ODC}=\widehat{OCD}=45^{\circ}$.\\

The diagonals of square are equal, so half of the diagonals are equal too. From this statement, $u=w$, $o=p_{1}$, $q_{1}=d_{2}$, $e_{2}=q$. Then, we can write the
isosceles triangles $\Delta ODC=\Delta PDB=\Delta BQA=\Delta AH_{1}D=\Delta DOC$. Since, square has equal sides, then $g=f=e=d$. Using the special triangle (45-45-90), $\frac{g\sqrt{2}}{2}=w$, $\frac{f\sqrt{2}}{2}=o$, $\frac{e\sqrt{2}}{2}=d_{2}$, $\frac{d\sqrt{2}}{2}=q$ and $\frac{g\sqrt{2}}{2}=w=o=d_{2}=q$.
Thus, $|OP|=|PQ|=|QH_{1}|=|H_{1}O|$, as we desired.\\

The ratio between the lengths of edges of the main square and the discovered/new square:\\

Let's think that one side of the main square is $x$, $(x\in \mathbb{Z})$.
\begin{eqnarray*}
&&\frac{\sin90^{\circ}}{x}=\frac{\sin45^{\circ}}{|OC|}\\
&&\frac{\sqrt{2}}{2}x=|OC|.
\end{eqnarray*}

Since $|OP|$ is the double of $|DC|$ because of edges $w=o$, $|OP|=\sqrt{2}x$.
The one side of the main square is $x$, then one side of the discovered square
need to be $\sqrt{2}x$ since the ratio between the sides will not change no matter the size of the square is because same sine law will be applied with the same angles.

\subsection{Application to Octagon}

\begin{figure}

  \includegraphics[width=120mm]{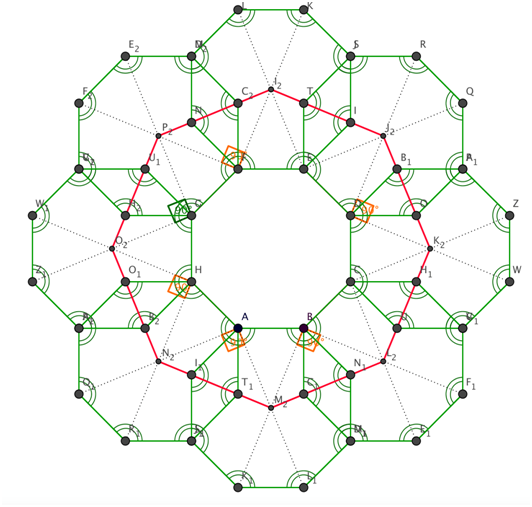}\\
  \caption{}\label{Fig. 5}

\end{figure}

The main octagon is $ABCDEFGH$ and the eight octagons, which are outside the main octagon(Figure \ref{Fig. 5}).\\

Now, we know that, the sum of the interior angles of any octagon is $1080^{\circ}$ and the diagonals are bisectors of the octagon, which are $67,5^{\circ}$.
 Thus, we obtain that each of interior angle of a regular octagon is $135^{\circ}$. And also, we can say that
$\widehat{I_{2}FE}=\widehat{I_{2}EF}=<\widehat{J_{2}ED}=<\widehat{J_{2}DE}=<\widehat{K_{2}CD}=<\widehat{K_{2}DC}=<\widehat{L_{2}CB}=<\widehat{L_{2}BC}=<\widehat{BM_{2}A}=<\widehat{AM_{2}B}
=<\widehat{AN_{2}H}=<\widehat{HN_{2}A}=<\widehat{HO_{2}C}=<\widehat{CO_{2}H}=<\widehat{CP_{2}F}=<\widehat{FP_{2}C}$.\\

Diagonals of regular octagon are equal. So half of the diagonals are equal too. Therefore, $|I_{2}F|=|I_{2}E|=|J_{2}E|=|J_{2}D|=|K_{2}D|=|K_{2}C|=|L_{2}C|=|L_{2}B|=|BM_{2}|=|AM_{2}|=|N_{2}A|=|N_{2}H|
=|O_{2}H|=|O_{2}C|=|P_{2}C|=|P_{2}F|$. Since they see the same angle, which is $45^{\circ}$, then $|CF|=|FE|=|ED|=|DC|=|CB|=|BA|=|AH|=|HC|$. The triangles $\Delta AM_{2}B=\Delta BL_{2}C=\Delta CK_{2}D=\Delta DJ_{2}F=\Delta FI_{2}E=\Delta CP_{2}F=\Delta HO_{2}C=\Delta AN_{2}H$ are isosceles triangles. This step completes the proof.\\

The ratio between the lengths of edges of the main octagon and the discovered/new octagon:\\

Assume that one side of the main regular octagon is $x$, $(x\in \mathbb{Z})$.
\begin{eqnarray*}
&&\frac{\sin45^{\circ}}{x}=\frac{\sin67,5^{\circ}}{|K_{2}C|}\\
&&1,31x=|K_{2}C|.
\end{eqnarray*}
Also $|K_{2}C|$ and $|CL_{2}|$ are equal because they see the same angle which is $67,5^{\circ}$. If we say the edges which see $45^{\circ}$ will be $x$ then the edges which see $90^{\circ}$ will be $\sqrt{2}x$.
\begin{eqnarray*}
&&|K_{2}C|=|CL_{2}|=1,31x\\
&&|K_{2}L_{2}|=1,85x.
\end{eqnarray*}
The one side of the main octagon is $x$ then, one side of the discovered octagon need to be $1.85x$ since the ratio between sides won't change no matter the size of the hexagon because same sine law will be applied with the same, equal angles.

\section{Conclusion}

In this work, we have investigated how Napoleon's Theorem is applied to regular polygons.
Firstly, we gave a simple proof of the Theorem using Sine and Cosine Laws.
Further, we applied the Theorem to regular polygons such as hexagon, square and octagon, with the same idea.
We draw a triangle, which is not equilateral, and then draw equilateral triangles outside the triangle. We got
the centroids of the equilateral triangle, unit them with lines.\\

In our application on hexagon, we found that the main
regular hexagon was  $x$ and the discovered square was $\sqrt{3}x$. Hence, the main regular hexagon would every time
be nearly twice as small as the new regular hexagon. We did the same ratios in the application of square. This times the main square
was $x$, and the discovered square was $\sqrt{2}x$. Hence, this time the main square was one and a half smaller than the new square.
And also, we used the octagon. Thus, the regular polygons length ratios which are the main and
discovered polygon, was transferred in graphs and seen that there was a strong correlation between them.

\end{document}